\documentclass[12pt]{article}
\usepackage{amsmath, amssymb}
\usepackage[pdftex]{graphicx}
\usepackage{psfrag,epsf}
\usepackage{url} % not crucial - just used below for the URL 
\usepackage{natbib} % no se que hace
\usepackage[pagebackref, colorlinks=true, urlcolor=blue, citecolor=blue, linkcolor=blue]{hyperref}

\usepackage{setspace}
\def\I{{\mathbb I}}
\def\I{{\mathbb I}}
\def\P{{\mathbb P}}
\def\E{{\mathbb E}}

\def\y{{\vec{y}}}
\def\Y{{\vec{Y}}}

\usepackage{mathtools}
\newcommand\givenbase[1][]{\! \:#1\lvert\: \!}
\let\given\givenbase

\DeclarePairedDelimiterX\Basics[1](){\let\given\sgiven #1}

\newtheorem{thm}{Theorem}

\newtheorem{Cor}{Corollary}
\newtheorem{Lemma}{Lemma}
\newtheorem{Assumption}{Assumption}
\newtheorem{Example}{Example}

\newtheorem{Notation}{Notation}
\newtheorem{Definition}{Definition}

\usepackage{changes}
\setremarkmarkup{(#2)}

\renewcommand{\vec}[1]{\mathbf{#1}}

\title{Teaching decision theory proof strategies
       using a crowdsourcing problem
			}
\author{
	Lu\'is Gustavo Esteves							\\
	University of S\~ao Paulo					\\ 
	\and   Rafael Izbicki							\\
	Federal University of S\~ao Carlos \\
	\and Rafael Bassi Stern						\\
	Federal University of S\~ao Carlos \\
}
\date{\today}

\begin{document}

\def\spacingset#1{\renewcommand{\baselinestretch}%
{#1}\small\normalsize} \spacingset{1}

\maketitle

\begin{abstract}
Teaching how to derive minimax decision rules can be challenging because of the lack
of examples that are simple enough to be used in the classroom. Motivated by this
challenge, we provide a new example  that illustrates
the use of standard techniques in
the derivation of optimal decision rules
under the Bayes and minimax approaches.
We discuss how to predict the value of an unknown quantity, $\theta \! \in \! \{0,1\}$,
given the opinions of $n$ experts.
An important example of such crowdsourcing problem occurs in modern
cosmology, where $\theta$ indicates whether a given galaxy is merging or not,
and $Y_1, \ldots, Y_n$ are the opinions 
from $n$ astronomers regarding $\theta$.
We use the obtained prediction rules 
to discuss advantages and disadvantages of
the Bayes and minimax approaches to decision theory.
The material presented here is intended to be taught to 
first-year graduate students.
\end{abstract}


\


\

{\bf Running title: Decision theory using a crowdsourcing problem}  

\newpage
\spacingset{1.45} % DON'T change the spacing!

\section{Introduction}

Decision theory formally compares statistical methods. 
It is deeply linked to most Bayesian procedures \citep{casella2002statistical,degroot2005optimal,parmigiani2009decision}.
Also, in the minimax approach, decision theory is 
the standard criterion for defining frequentist optimality
of nonparametric methods \citep{wasserman2006all,tsybakov2008introduction}.
However, there are few published examples
in which both the minimax and Bayes solutions 
involve elementary calculations, 
which tends to obscure teaching this topic.
Here, we explore a motivating example where
calculations are simple and yet illustrate 
standard techniques for finding optimal estimators.
We also use this example to discuss 
advantages and disadvantages of 
the Bayes and minimax approaches to decision theory.
The material presented here is intended to be used for first-year graduate level students.

The main element of the motivating example is a sequence, 
$\Y = \{Y_1, \ldots, Y_n\} \! \in \! \{0,1\}^{n}$,
composed of the opinions of $n$ experts about an unknown quantity, $\theta \in \{0,1\}$.
For example, $\theta$ might be 1 if a given galaxy is merging with another galaxy and $0$ otherwise.
It is an important question in cosmology to infer the value of $\theta$, 
because the number of merger galaxies plays a key role in testing theories about
the evolution of the Universe \citep{lotz2004new,freeman2013new}.
In this case, $Y_1, \ldots, Y_n$ would be the opinions 
from $n$ astronomers regarding the value of $\theta$, obtained using an image of that galaxy;
see Table \ref{tbl:galaxy} for
some examples from the Cosmic Assembly Near-IR Deep Extragalactic Legacy
Survey (CANDELS; \citealt{0067-0049-197-2-35}). 
Given the low resolution of the images, it is common for astronomers to disagree regarding the value of $\theta$. 
Hence, we are interested in using all the information in the opinions, $y_1, \ldots, y_n$, 
to infer the value of $\theta$. 
This example is a building block for crowdsourcing models,
such as the one developed in \cite{izbicki2013learning},
which are becoming increasingly popular due to platforms such as Amazon Mechanical Turk
and projects such as Galaxy Zoo\footnote{\url{https://www.galaxyzoo.org/}} and foldit\footnote{\url{http://fold.it}}.

As a minor element of the motivating example,
we also allow the use of a fair coin flip, $U$, 
to infer about $\theta$.
That is, $U \sim \text{Bernoulli}(0.5)$ and $U$ is independent of $\Y$.
The coin flip can be useful in making a decision when 
one is indifferent about the options that are available.
The importance of this element is discussed in more detail in subsection \ref{sec::minimaxFixed}.
%We study the problem
%in two settings: the first one is to consider $\theta$ to be a fixed but unknown quantity, that is, a parameter in the traditional frequentist sense. The second approach is to consider $\theta$ to be a random variable. To maintain agreement with the traditional frequentist notation, in the latter case instead of using the character $\theta$, we use $Z$ to represent this quantity.\footnote{ 
%	This is often referred to  as  a latent variable \red{or just random variable?}, that is, an unknown quantity to which we assign probability measure.}

\begin{table}[h!]
  \centering
  \begin{tabular}{ | c |  c  c  c | }
    \hline
  &  \multicolumn{3}{c|}{Is the galaxy a merger?} \\ \hline
    Galaxy & Expert 1 & Expert 2  & Expert 3 \\ \hline
    \begin{minipage}{.1\textwidth}
    \center
      \includegraphics[width=\linewidth, height=18mm,width=18mm,page=1,trim={10cm 10cm 0.5cm 1cm},clip]{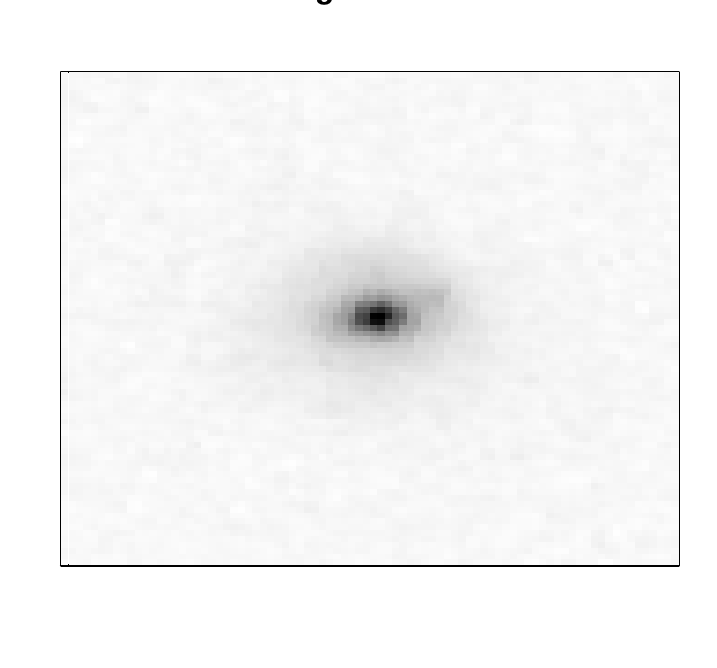}
    \end{minipage}
    &
    No
    %\begin{minipage}[t]{5cm}
    %\end{minipage}
    & 
    Yes
    & 
    Yes
    %\begin{minipage}{5cm}
    %\end{minipage}
    \\ \hline
    \begin{minipage}{.1\textwidth}
    \center
      \includegraphics[width=\linewidth, height=18mm,width=18mm,page=11,trim={10cm 10cm 0.5cm 1cm},clip]{exampleFalseMerAllBandsRealSizeWithWithoutMask.pdf}
    \end{minipage}
    &
    No
    %\begin{minipage}[t]{5cm}
    %\end{minipage}
    & 
    Yes
    & 
    No
    %\begin{minipage}{5cm}
    %\end{minipage}
    \\ \hline
    \begin{minipage}{.1\textwidth}
    \center
\includegraphics[width=\linewidth, height=18mm,width=18mm,page=37,trim={10cm 10cm 0.5cm 1cm},clip]{exampleFalseMerAllBandsRealSizeWithWithoutMask.pdf}
    \end{minipage}
    &
    Yes
    %\begin{minipage}[t]{5cm}
    %\end{minipage}
    & 
    Yes
    & 
    No
    %\begin{minipage}{5cm}
    %\end{minipage}
    \\ \hline
  \end{tabular}
  \caption{An example of crowdsourcing in cosmology: three astronomers give their opinion about weather a given galaxy
  is a merger or not. What is the optimal way to combine such opinions in order to recover the truth? }\label{tbl:galaxy}
\end{table}

The motivating example can be framed as a decision problem,
hereafter called the crowdsourcing problem.
For each set of expert opinions and coin flip outcome,
the decision-maker must make a prediction, $a\in\{0,1\}$,
about the value of $\theta$.
Given values for $a$ and $\theta$,
the decision-maker is penalized according to 
$L(\theta,a)=\I(\theta\neq a)$, that is,
errors are penalized with a loss of one unit
and correct decisions have zero loss. 
The goal of the decision-maker is to design a decision rule
$\delta \! : \! \{0,1\}^n \times \{0,1\} \longrightarrow \{0,1\}$ which, 
for each collection of expert opinions $\y:=(y_1,\ldots,y_n)$ and coin flip outcome, $u$, 
associates a decision $\delta(\y, u) \in \{0,1\}$. 
The decision-maker wants a decision rule
that leads to a small loss.
More precisely, the risk \citep{wasserman2006all}[p.228] of the decision rule $\delta$,
$R(\delta,\theta)=\E_{\Y \given \theta}[L(\delta(\Y, U),\theta)]$,
should be small.
Because the risk depends on the unknown value of
$\theta$ and it is typically not possible to find a decision rule
that minimizes the risk for every $\theta$, 
the following optimality criteria are typically used:
\begin{itemize}
  \item \textbf{[Minimax]} $\delta^*$ is a minimax rule if 
	$\sup_{\theta} R(\delta^*,\theta)=\inf_{\delta}\sup_{\theta} R(\delta,\theta)$
  \item \textbf{[Bayes]} $\delta^*$ is a Bayes decision rule if 
	$\E_{P}[R(\delta^*,\theta)]=\inf_{\delta}\E_{P}[ R(\delta,\theta)]$
\end{itemize}

$\E_{P}[ R(\delta,\theta)]$ is called the Bayes risk of $\delta$.
The expectation that appears in the Bayes risk is with
respect to $P$, the prior distribution on $\theta$.

In order to explore the crowdsourcing problem, 
we require some assumptions:
\begin{Assumption} \
 \label{assumption:accuracy}
 \begin{enumerate}
  \item $Y_{1},\ldots,Y_{n}$ are independent (conditionally on $\theta$).
  \item $\gamma_i := \P(Y_i=1 \given  \theta=1)=\P(Y_i=0 \given   \theta=0)$.
	\item $n$ is odd.
 \end{enumerate}
\end{Assumption}

Observe that $\gamma_i$ is
the accuracy of expert $i$.
According to Assumption \ref{assumption:accuracy}.2, for each expert, 
the chance of making a mistake does not depend 
on the value of $\theta$.
Assumption \ref{assumption:accuracy}.3 simplifies some derivations
by considering that one cannot obtain ties between experts' 0's and 1's.

Section \ref{sec::known} analyzes the case in which $\gamma_i$'s are known in advance.
Section \ref{sec::unknown} analyzes the case in which the $\gamma_i$'s are unknown.
In pursuance of simplifying the text, these sections use the following notation:

\begin{Notation} \
 \begin{enumerate}
  \item $c := P(\theta=1)$.
  \item $\bar{\Y} := \frac{1}{n}\sum_{i=1}^{n} Y_i$ and $\bar{\y} := \frac{1}{n}\sum_{i=1}^{n}y_i$.
 \end{enumerate}
\end{Notation}

\section{Experts' accuracies are known}
\label{sec::known}

\subsection{Bayes decision rule}
\label{sec::bayesFixed}

 This subsection uses the crowdsourcing problem with known expert accuracies 
 to illustrate a strategy that is commonly used to find Bayes decision rules.
 In order to simplify calculations and make the Bayes decision rules more interpretable, 
 we use the following notation:
\begin{Notation}
 \label{notation:known}
  $w_i^* := 2\log\left( \frac{ \gamma_i}{1-\gamma_i}\right)$ and
 $C^* := \log\left(\frac{1-c}{c}\right)+ \sum_{i=1}^{n} \frac{w_i^{*}}{2}$.
\end{Notation}
Observe that $w_i^*$ increases monotonically with the $i$-th expert's accuracy, $\gamma_{i}$.
Furthermore, if $w_i^* < 0$, then $\gamma_{i} < \frac{1}{2}$, that is, 
the $i$-th expert is less accurate than a coin flip.
Similarly, if $w_i^* > 0$, then the $i$-th expert is more accurate than a coin flip.
Next, we derive Bayes rules in the crowdsourcing problem.

\begin{thm}
 \label{theorem:weighted_experts}
  In the crowdsourcing problem with known expert accuracies,
  for each $c \in [0,1]$, the following decision rule is a Bayes decision rule
	\begin{align*}
		\delta^{*}_c(\y,u)	&=
		\begin{cases}
			1 & \text{, if } \sum_{i=1}^{n} w_i^* y_i > C^*	\\
			0 & \text{, if } \sum_{i=1}^{n} w_i^* y_i < C^*	\\
			u	& \text{, if } \sum_{i=1}^{n} w_i^* y_i = C^*
		\end{cases}
	\end{align*}
\end{thm}
The decision rule given by Theorem \ref{theorem:weighted_experts} admits an intuitive interpretation. 
If a weighted sum of the opinions of the experts is greater than the cutoff, $C^{*}$,
then one decides that $\theta$ is 1.
$C^*$ varies according to the prior probability assigned to $\theta$. 
Also, the more an expert is accurate, 
the larger the weight assigned to him.
Observe that, if expert $i$ is worse than a coin flip ($w_i^* < 0$),
then the model favors $\theta=1-y_{i}$, the opposite of the expert's opinion.
Roughly speaking, the decision will be determined by the most accurate experts.

In order to prove Theorem \ref{theorem:weighted_experts},
a commonly used technique is given by Theorem \ref{theorem:extensive_form}:
\begin{thm}[\citep{parmigiani2009decision}]
 \label{theorem:extensive_form}
 A decision rule $\delta^{*}$ is a Bayes decision rule
 if and only if, for each $\y$ and $u$, 
 $\delta^{*}(\y, u)$ minimizes
 $\E_{P}[L(\delta(\y, u),\theta)|\y]$.
\end{thm}

Theorem \ref{theorem:extensive_form} allows one
to find a Bayes decision rule by determining a decision
for one value of $\y$ at a time.
This procedure is called the extensive form of a decision problem.
It is generally easier to solve a problem in the extensive form
than to directly prove that 
a given decision rule has a Bayes risk that is 
smaller than every other decision rule  (the latter approach is called the normal form of a decision problem).
The following paragraphs show how the solution of the extensive form 
for the crowdsourcing problem can be simplified to Theorem \ref{theorem:weighted_experts}.

In the crowdsourcing problem $\delta(\textbf{y}, u) \in \{0,1\}$.
Using this observation, one can obtain
from Theorem \ref{theorem:extensive_form}
the following lemma:
\begin{Lemma}
 \label{Lemma:1}
 $\delta^{*}$ is a Bayes decision rule in the crowdsourcing problem if and only if, for all $\y$ and $u$,
 \begin{align*}
	\begin{cases}
	 \delta^{*}(\y, u) =1 & \text{, if } \frac{\E_P[L(0,\theta) \given \y]}{\E_P[L(1,\theta) \given \y]} > 1	\\
	 \delta^{*}(\y, u) =0 & \text{, if } \frac{\E_P[L(0,\theta) \given \y]}{\E_P[L(1,\theta) \given \y]} < 1
	\end{cases}
 \end{align*}
\end{Lemma}

One can conclude from Lemma \ref{Lemma:1} that 
there may exist several Bayes decision rules for the crowdsourcing problem.
If there exists some $\textbf{y}$ such that
$\frac{\E_P[L(0,\theta) \given \y]}{\E_P[L(1,\theta) \given \y]} = 1$,
then one is indifferent between choosing either $0$ or $1$
when observing $\textbf{y}$.
In other words, in this case 
the choice that is made is not relevant for the goal of obtaining a Bayes decision rule.
Corollary \ref{Cor:bayes_flip}
presents a particular Bayes decision rule
that aligns with the one in Theorem \ref{theorem:weighted_experts}.
This case of Bayes decision rule is emphasized because
it is used in subsection \ref{sec::minimaxFixed}.
\begin{Cor}
 \label{Cor:bayes_flip}
 Let $\delta^{*}$ be defined as
 \begin{align*}
	\delta^{*}(\y, u) =
	\begin{cases}
	 1 & \text{, if } \frac{\E_P[L(0,\theta) \given \y]}{\E_P[L(1,\theta) \given \y]} > 1	\\
	 0 & \text{, if } \frac{\E_P[L(0,\theta) \given \y]}{\E_P[L(1,\theta) \given \y]} < 1	\\
	 u & \text{, if } \frac{\E_P[L(0,\theta) \given \y]}{\E_P[L(1,\theta) \given \y]} = 1
	\end{cases}
 \end{align*}
 $\delta^{*}$ is a Bayes decision rule for the crowdsourcing problem.
\end{Cor}

Also observe that  $L(0,0)=L(1,1)=0$ and $L(0,1)=L(1,0)=1$.
Therefore, one can obtain Lemma \ref{Lemma:2}
\begin{Lemma}
 \label{Lemma:2}
 \begin{align*}
  \begin{cases}
	 \frac{\E_P[L(0,\theta) \given \y]}{\E_P[L(1,\theta) \given \y]} > 1 &\iff \frac{\P(\theta=1 \given \y)}{\P(\theta=0 \given \y)} > 1 \iff	\P(\theta=1 \given \y) > \frac{1}{2} \\
	 \frac{\E_P[L(0,\theta) \given \y]}{\E_P[L(1,\theta) \given \y]} < 1 &\iff \frac{\P(\theta=1 \given \y)}{\P(\theta=0 \given \y)} < 1 \iff \P(\theta=1 \given \y) < \frac{1}{2}
	\end{cases}
 \end{align*}
\end{Lemma} 

Furthermore, 
\begin{align}
 \label{equation:1}
 \P(\theta=1 \given \y)>\frac{1}{2} &\iff \frac{c \prod_{i=1}^{n} (\gamma_i)^{y_i} (1-\gamma_i)^{1-y_i} }{c \prod_{i=1}^{n} (\gamma_i)^{y_i} (1-\gamma_i)^{1-y_i}+(1-c)\prod_{i=1}^{n} (1-\gamma_i)^{y_i} (\gamma_i)^{1-y_i}} >\frac{1}{2} \nonumber \\
%&\iff \frac{c}{1-c}\frac{\prod_{i=1}^{n}(\gamma_i)^{y_i} (1-\gamma_i)^{1-y_i}}{\prod_{i=1}^{n}(1-\	gamma_i)^{y_i} (\gamma_i)^{1-y_i}} >  1 \\
 &\iff  2\sum_{i=1}^{n} \log\left( \frac{ \gamma_i}{1-\gamma_i}\right) y_i >\log\left(\frac{1-c}{c}\right)+ \sum_{i=1}^{n} \log\left(\frac{\gamma_i}{1-\gamma_i}\right).
\end{align}
%
%
%$$\log\left(\frac{c}{1-c}\right)+\log\left[ \prod_{i=1}^{n}\left(\frac{\pi_i}{1-\gamma_i}\right)^{y_i} \left(\frac{1-\pi_i}{\gamma_i}\right)^{1-y_i}\right]>0 \iff$$
%
%$$\log\left(\frac{c}{1-c}\right)+\sum_i \log\left[ \left(\frac{\pi_i \gamma_i}{(1-\gamma_i)(1-\pi_i)}\right)^{y_i} \left(\frac{1-\pi_i}{\gamma_i}\right)\right]>0 \iff  $$ 

Plugging Notation \ref{notation:known} into Equation \ref{equation:1}, one obtains that
\begin{align*}
	\P(\theta=1 \given \y)>\frac{1}{2} &\iff \sum_{i=1}^{n} w_i^* y_i > C^*
\end{align*}

Similarly,
\begin{align*}
	\P(\theta=1 \given \y)<\frac{1}{2} &\iff \sum_{i=1}^{n} w_i^* y_i < C^*
\end{align*}

Combining the above derivation with
Corollary \ref{Cor:bayes_flip} and Lemma \ref{Lemma:2},
Theorem \ref{theorem:weighted_experts} is obtained.

Besides the importance of Bayes decision rules from the Bayesian perspective,
these decision rules can also have good frequentist properties. This is exemplified in the next subsection, which
illustrates how to find a minimax decision rule 
through Bayes decision rules. 
This proof strategy is commonly used.

\subsection{Minimax decision rule}
\label{sec::minimaxFixed}

This subsection exemplifies a strategy that is often used to
find minimax decision rules. In particular, it proves 
the following theorem in the crowdsourcing problem:
\begin{thm}
 \label{theorem:minimax_fixed}
 $\delta^{*}_{\frac{1}{2}}$ (theorem \ref{theorem:weighted_experts}) 
 has constant risk and is minimax.
\end{thm}
Observe that, if $\P(\theta=1) = c = \frac{1}{2}$,
then one is indifferent between $\theta=1$ and $\theta=0$
before observing the experts' opinions.
This is a possible intuition for why it is plausible that
the Bayes decision rule $\delta^{*}_{\frac{1}{2}}$ is also minimax.
Furthermore, Corollary \ref{corollary:minimax_majority} shows how 
$\delta^{*}_{\frac{1}{2}}$ can be simplified into a familiar decision rule.
\begin{Cor}
 \label{corollary:minimax_majority}
 If $\gamma_1 = \ldots = \gamma_n>\frac{1}{2}$,
 that is, all experts have the same accuracy and
 the accuracy is better than that of a coin flip, 
 then the decision rule $\delta_{M}$ defined as
 \begin{align*}
   \delta_{M}(\y, u) = 1 \iff \bar{\y} > \frac{1}{2}
 \end{align*}
 is equal to the Bayes decision rule $\delta^{*}_{\frac{1}{2}}(\y, u)$ and
 is a minimax decision rule.
 Observe that $\bar{\y}=\frac{1}{2}$ cannot occur because 
 of Assumption \ref{assumption:accuracy}.3.
\end{Cor}
The decision rule $\delta_{M}$ is known 
as the majority rule and has been used in crowdsourcing problems \citep{raykar2010learning}.
Such an agreement between $\delta^{*}_{\frac{1}{2}}$ and $\delta_{M}$
might also improve the plausibility of $\delta^{*}_{\frac{1}{2}}$ being minimax.

In order to prove that $\delta^{*}_{\frac{1}{2}}$ is indeed minimax,
we use a standard technique, that is described in Theorem \ref{theorem:minimax_technique}.
\begin{thm}[\citep{wasserman2006all}]
  \label{theorem:minimax_technique}
  If $\delta^{*}$ is a Bayes decision rule against some prior distribution and
  the risk of $\delta^{*}$, $R(\delta^{*},\theta)$,
  is constant on $\theta$, then $\delta^{*}$ is a minimax decision rule.
\end{thm}

Theorem \ref{theorem:minimax_technique}
provides a proof technique for finding minimax rules that 
is usually better than directly applying the definition of minimax decision rule.
Observe that a minimax decision rule attains an infimum risk
over the space of all possible decision rules.
This space can be large and, as a result,
it is generally complicated to directly perform
calculations for all the elements in it. 
Theorem \ref{theorem:minimax_technique} allows one
to perform calculations for a single decision rule and 
prove that this ``candidate" is minimax. 

In the crowdsourcing problem,
Theorem \ref{theorem:weighted_experts}
obtains a Bayes decision rule $\delta^{*}_{c}$
for each prior $c = P(\theta=1)$.
Hence, if some $\delta^{*}_{c}$
has constant risk, it follows from Theorem \ref{theorem:minimax_technique}
that $\delta^{*}_{c}$ is minimax. Observe that

\begin{align*}
  \begin{cases}
   R(\delta^{*}_c,0) &= \P\left(\sum_{i=1}^{n} w_i^* Y_i > C^* \given[\Big]  \theta=0\right) + \frac{1}{2}\P\left(\sum_{i=1}^{n} w_i^* Y_i = C^* \given[\Big]  \theta=0\right)	\\
   R(\delta^{*}_c,1) &= \P\left(\sum_{i=1}^{n} w_i^* Y_i < C^* \given[\Big]  \theta=1\right)	+ \frac{1}{2}\P\left(\sum_{i=1}^{n} w_i^* Y_i = C^* \given[\Big]  \theta=1\right)	
  \end{cases}
\end{align*}

Hence, by combining Theorems \ref{theorem:weighted_experts} and \ref{theorem:minimax_technique}
one obtains 

\begin{Lemma}
 \label{Lemma:3}
 If 
 \begin{align*}
  &\P\left(\sum_i w_i^* Y_i > C^* \given[\Big]  \theta=0\right) + \frac{1}{2}\P\left(\sum_{i=1}^{n} w_i^* Y_i = C^* \given[\Big]  \theta=0\right)	\\
	=&\P\left(\sum_{i=1}^{n} w_i^* Y_i < C^* \given[\Big]  \theta=1\right)	+ \frac{1}{2}\P\left(\sum_{i=1}^{n} w_i^* Y_i = C^* \given[\Big]  \theta=1\right)	
 \end{align*}
 then $\delta^{*}_{c}$ has constant risk and is minimax.
\end{Lemma}

In order to find a $\delta^{*}_{c}$ such as in Lemma \ref{Lemma:3},
one must choose an appropriate value for $c$.
It follows from Assumption \ref{assumption:accuracy} that
$\Y|\theta=0$ is identically distributed to $1-\Y|\theta=1$.
This fact is used from the second to the third line of the following derivation:

\begin{align}
\label{eq::equalityRisk}
&\P\left(\sum_{i=1}^{n} w_i^* Y_i > C^* \given[\Big]  \theta=0\right) + \frac{1}{2}\P\left(\sum_{i=1}^{n} w_i^* Y_i = C^* \given[\Big]  \theta=0\right)	\notag\\
=& \P\left(\sum_{i=1}^{n} w_i^* \left(Y_i-\frac{1}{2}\right) > \log\left(\frac{c}{1-c}\right) \given[\Big]  \theta=0\right) + \frac{1}{2}\P\left(\sum_{i=1}^{n} w_i^* \left(Y_i-\frac{1}{2}\right) = \log\left(\frac{c}{1-c}\right) \given[\Big]  \theta=0\right)	\notag\\
=& \P\left(\sum_{i=1}^{n} w_i^* \left(\frac{1}{2}-Y_i\right) > \log\left(\frac{c}{1-c}\right)  \given[\Big]  \theta=1\right) + \frac{1}{2}\P\left(\sum_{i=1}^{n} w_i^* \left(\frac{1}{2}-Y_i\right) = \log\left(\frac{c}{1-c}\right)  \given[\Big]  \theta=1\right)	\notag\\
=& \P\left(\sum_{i=1}^{n} w_i^* \left(Y_i-\frac{1}{2}\right) < -\log\left(\frac{c}{1-c}\right)  \given[\Big]  \theta=1\right) + \frac{1}{2}\P\left(\sum_{i=1}^{n} w_i^* \left(Y_i-\frac{1}{2}\right) = -\log\left(\frac{c}{1-c}\right)  \given[\Big]  \theta=1\right)	\notag\\
\end{align}
Moreover,
\begin{align}
\label{eq::equalityRisk_2}
&\P\left(\sum_{i=1}^{n} w_i^* Y_i < C^* \given[\Big]  \theta=1\right)	+ \frac{1}{2}\P\left(\sum_{i=1}^{n} w_i^* Y_i = C^* \given[\Big]  \theta=1\right)\notag\\
=& \P\left(\sum_{i=1}^{n} w_i^* \left(Y_i-\frac{1}{2}\right) < \log\left(\frac{c}{1-c}\right)  \given[\Big]  \theta=1\right) + \frac{1}{2}\P\left(\sum_{i=1}^{n} w_i^* \left(Y_i-\frac{1}{2}\right) = \log\left(\frac{c}{1-c}\right)  \given[\Big]  \theta=1\right)	\notag\\	
\end{align}

Observe that the last term in Equation \ref{eq::equalityRisk}
equals to the last term in Equation \ref{eq::equalityRisk_2} when
$\log(\frac{c}{1-c}) = -\log(\frac{c}{1-c}) = 0$.
This occurs when $c=\frac{1}{2}$.
Theorem \ref{theorem:minimax_fixed} is proved from combining Lemma \ref{Lemma:3} and 
Equations \ref{eq::equalityRisk} and \ref{eq::equalityRisk_2}.
The coin flip, $U$, was essential in the derivation
of the minimax decision rule,
since it allows one to divide the risk equally between $\theta=1$ and $\theta=0$
when $\sum_{i=1}^{n} w_i^* Y_i = C^*$.

\subsection{Least favourable prior}

The crowdsourcing problem can also be used to illustrate
a strategy for finding a \emph{least favorable prior}.

\begin{Definition}
 \label{definition:least_favorable}
 Let $P^*$ be a prior distribution associated with a Bayes rule, $\delta^{*}$.
 $P^*$ is said to be a least favorable prior if, for every other prior, $P$,
 and associated Bayes rule, $\delta$,
 \begin{align*}
  E_{P^{*}}[R(\delta^{*}, \theta)] \geq E_{P}[R(\delta, \theta)]
 \end{align*}
\end{Definition}
In words, the least favorable prior
leads one to use a Bayes rule with the largest Bayes risk.
In this sense, the least favorable prior
is such that it makes it the hardest to obtain a desirable result
in the decision problem.
In order to find a least favorable prior, 
one can use the following theorem:
\begin{thm}[\citep{parmigiani2009decision}]
 \label{theorem:least_favorable}
 A prior associated with a Bayes rule with constant risk
 is a least favorable prior.
\end{thm}
 
One can combine Theorems \ref{theorem:minimax_fixed} and \ref{theorem:least_favorable}
to conclude that $P^*(\theta=1)=1/2$
is the least favorable prior in the crowdsourcing problem.
Such a use of Theorem \ref{theorem:least_favorable} is usually
easier than directly exploring the full space of prior distributions
for a least favorable prior. 
%Indeed, it is the prior that reflects less knowledge about 
%$\theta$ in Laplace's sense.

%Notice however that this is not necessarily the only prior that leads to a Bayes decision rule with constant risk in our problem. 
%This is because the experts'  guesses $Y_i$ assume a discrete value, and therefore one can change $c$
%and still get the same decision rule as $\delta_{1/2}$. Notice, however, that the posterior probabilities 
%$\P(\theta|\y)$ do change for each $c$. Finally, because
%$\delta_{1/2}$ is the only Bayes rule with respect to
%$\P*$, $\delta_{1/2}$ is the only minimax rule \citep{parmigiani2009decision}.

%Hence, the majority vote is a minimax decision in this setting. 
%\red{mencionar? Note that Condorcet proved that this decision dominates the one of letting
%a single expert decides if $\theta=1$ or $\theta=0.$}

\section{Experts' accuracies are unknown}
\label{sec::unknown}

Although the previous section considers the experts' accuracies to be known,
one often is uncertain about these quantities.
In the following subsections, 
we consider the case in which the experts' accuracies are unknown parameters
and use this case to illustrate more proof strategies.
Formally, in subsections \ref{section:bayes_unknown} and \ref{section:minimax_unknown} 
the parameter space is
\begin{align*}
  \Theta=\{(\theta,\gamma_1, \ldots, \gamma_n):\mbox{ }\theta \in \{0,1\}, \gamma_i \in [0,1] \mbox{ }\forall i=1,\ldots,n\}
\end{align*}

\subsection{Bayes decision rule}
\label{section:bayes_unknown}

%Any decision rule is minimax in this setting. 
%To note this, first note that the risk of this
%decision rule is $\frac{1}{2}$ for both values of $\theta$. Now, take any other not randomized decision.
%Note that if $R(\delta,0) =P(\delta(Y)=0 \given[\Big]  \theta=1,\pi_1, \ldots, \pi_n,\gamma_1, \ldots, \gamma_n)<\frac{1}{2}$ ,we
%must have $R(\delta,0) =P(\delta(Y)=0 \given[\Big]  \theta=1,\pi_1, \ldots, \pi_n,\gamma_1, \ldots, \gamma_n)<\frac{1}{2}$

Considering the generalization of the parameter space,
not every prior allows an analytical derivation of the posterior distribution.
In order to be able to derive the Bayes decision rules in this case,
we consider only priors of the following form:

\begin{Assumption} \
 \label{assumption:bayes_unknown}
 \begin{enumerate}
  \item $\theta, \gamma_1, \ldots, \gamma_n$ are independent.
  \item $\gamma_i \sim \text{Beta}(\alpha_i, \beta_i)$ ($\alpha_i$ and $\beta_i$ are hyperparameters chosen \emph{a priori}).
 \end{enumerate}
\end{Assumption} 

Assumption \ref{assumption:bayes_unknown}.1 breaks the derivation of the posterior into smaller parts.
Assumption \ref{assumption:bayes_unknown}.2 yields a conjugate distribution for Bernoulli trials,
such as the crowdsourcing problem.

Under Assumption \ref{assumption:bayes_unknown},
the Bayes decision rules derived under unknown expert accuracies
can be compared to the ones derived when the accuracies are known.
Notation \ref{notation:unknown} is analogous to notation \ref{notation:known}
and is useful for concisely presenting these decision rules.

\begin{Notation}
 \label{notation:unknown}
 $w_i^{**}	:= 2\log\left(\frac{\alpha_i}{\beta_i}\right)$ and
 $C^{**}	:= \log\left(\frac{1-c}{c}\right) + \sum_{i=1}^{n}\frac{w_i^{**}}{2}$.
\end{Notation}

\begin{thm}
 \label{theorem:weighted_experts_unknown}
  For each $c \in [0,1]$ the following decision rule is a Bayes decision rule 
	\begin{align*}
	 \delta^{**}_c(\y, u) &=
	 \begin{cases}
		1 & \text{, if } \sum_{i=1}^{n} w_i^{**} y_i > C^{**}	\\
		0 & \text{, if } \sum_{i=1}^{n} w_i^{**} y_i < C^{**}	\\
		u & \text{, otherwise}
	 \end{cases}
	\end{align*}
\end{thm}

Observe that $\delta^{**}_c$ is analogous to $\delta^*_c$, 
in Theorem \ref{theorem:weighted_experts}.
That is, if a weighted sum of the experts' opinions is greater than a given cutoff,
then $\delta^{**}_c$ decides for $1$.
Furthermore, the weights in each decision rule also admit similar interpretations.
It follows from Assumption \ref{assumption:bayes_unknown} that 
$E[\gamma_i] = \frac{\alpha_i}{\alpha_i + \beta_i}$ is the expected accuracy of expert $i$.
As a result, $\delta^{**}_c$ can be obtained from $\delta^{*}_c$
by substituting every instance of $\gamma_i$ by $E[\gamma_i]$.
Therefore, the larger the expected accuracy of expert $i$,
the more this expert is decisive in determining the value of $\delta^{**}_c$.

\begin{Example} 
 \normalfont
 \label{ex::priors}
 The Bayes decision rule that is obtained in Theorem \ref{theorem:weighted_experts_unknown}
 is affected by the prior distribution for the experts' accuracies, $\gamma_1,\gamma_2$ and $\gamma_3$. 
 For example, consider the following prior specifications:
 \begin{itemize}
  \item Prior 1: $\gamma_1 \sim \mbox{Beta}(5,2); \gamma_2,\gamma_3 \sim \mbox{Beta}(1,1)$. 
	               The decision-maker believes that the first expert is accurate, 
								 and he is indifferent about the others.
	\item Prior 2: $\gamma_1,\gamma_2 \sim \mbox{Beta}(5,2); \gamma_3 \sim \mbox{Beta}(1,1)$. 
	               The decision-maker believes that the first two experts are accurate, 
								 and he is indifferent about the third.
	\item Prior 3: $\gamma_1,\gamma_2 \sim \mbox{Beta}(5,2); \gamma_3 \sim \mbox{Beta}(2,5)$. 
	               The decision-maker believes that the first two experts are accurate, 
								 and that the third expert is a spammer.
	\end{itemize}	
	{\normalfont
		\begin{table}
			\begin{center}
			\label{tab::outcomes}
				\begin{tabular}{ c | c | c | c} \hline
					&   \multicolumn{3}{|c}{Bayes decision rule using}  \\
					Experts' opinions & Prior 1 & Prior 2 & Prior 3 \\
					\hline
					(0,0,0) & 0 & 0   & 0\\ 
					(1,0,0) & 1 & 0/1 & 1\\ 
					(0,1,0) & 0 & 0/1 & 1\\
					(0,0,1) & 0 & 0   &0 \\
					(1,1,0) & 1 & 1   & 1\\
					(1,0,1) & 1 & 0/1 &0 \\
					(0,1,1) & 0 & 0/1 &0 \\
					(1,1,1) & 1 & 1	  & 1\\
					\hline
				\end{tabular}
				\caption{Bayes decision rules that are obtained in Example \ref{ex::priors} for each of the priors that are considered. 
				        $0/1$ is used when both decisions are Bayes (i.e., they both have the same expected posterior loss).}
			\end{center}
		\end{table}
	}
	Furthermore, consider that, for every one of these priors,
	$c=\frac{1}{2}$.
	
  Table \ref{tab::outcomes} presents the Bayes decisions that are obtained 
	for each of the prior that were discussed above. 
  This table illustrates that, by choosing a prior that properly encodes one's uncertainty, 
  it is possible to obtain a decision rule that rationally combine the experts' opinions about $\theta$.
  For example, according to Prior $1$, expert $1$ is the only accurate one.
  As a result, the Bayes decision rule always agrees with expert $1$,
  whatsoever the other opinions are.
  Similarly, under Prior $2$ the only accurate experts are $1$ and $2$.
  In this case, the Bayes decision agrees with experts $1$ and $2$
  when they agree with each other.
  When experts $1$ and $2$ disagree, there exist Bayes decision rules
  that agree with either one of them.
  Finally, according to Prior $3$ experts $1$ and $2$ are accurate
  and expert $3$'s opinion is frequently incorrect.
  As a result, the Bayes decision rule agrees with experts $1$ and $2$
  when they agree, and decides (exactly) against expert $3$
  when experts $1$ and $2$ disagree.
\end{Example}

In order to prove Theorem \ref{theorem:weighted_experts_unknown},
one can solve the extensive form of the crowdsourcing problem, 
as presented in Theorem \ref{theorem:extensive_form}.
In this case, Lemma \ref{Lemma:2} can be applied.
Therefore, it is sufficient to derive the posterior distribution for $\theta$.

Let $ \boldsymbol{\gamma}=(\gamma_1,\ldots,\gamma_n)$.
The posterior distribution of $\theta$ is given by

\begin{align}
\P(\theta=1 \given \y)
&= \frac{\P(\theta=1, \Y=\y)}{\P(\Y=\y)}																		&						     \nonumber	\\
&= \frac{\P(\theta=1)}{\P(\Y=\y)}\int{\P(\Y=\y|\theta=1,\boldsymbol{\gamma})f(\boldsymbol{\gamma})d\boldsymbol{\gamma}}									& \text{Assumption \ref{assumption:bayes_unknown}.1} \nonumber 	\\
&= \frac{c}{\P(\Y=\y)}\int \prod_{i=1}^{n}(\gamma_i)^{y_i} (1-\gamma_i)^{1-y_i}												
\prod_{i=1}^{n}\frac{\Gamma(\alpha_i+\beta_i)}{\Gamma(\alpha_i)\Gamma(\beta_i)}(\gamma_i)^{\alpha_i-1} (1-\gamma_i)^{\beta_i-1} d\boldsymbol{\gamma}						&						     \nonumber	\\
&=  \frac{c}{\P(\Y=\y)}\prod_{i=1}^{n}\frac{\Gamma(\alpha_i+\beta_i)}{\Gamma(\alpha_i)\Gamma(\beta_i)} \int  (\gamma_i)^{\alpha_i+y_i-1} (1-\gamma_i)^{\beta_i-y_i}			
d\boldsymbol{\gamma}																					&						     \nonumber	\\
&= \frac{c}{\P(\Y=\y)}\prod_{i=1}^{n}\frac{\Gamma(\alpha_i+\beta_i)}{\Gamma(\alpha_i)\Gamma(\beta_i)}\frac{\Gamma(\alpha_i+y_i)\Gamma(\beta_i+1-y_i)}{\Gamma(\alpha_i+\beta_i+1)}	& \label{equation:4}						\\
\text{Similarly,} \notag \\
\P(\theta=0 \given \y)	&= 
\frac{1-c}{\P(\Y=\y)}\prod_{i=1}^{n}\frac{\Gamma(\alpha_i+\beta_i)}{\Gamma(\alpha_i)\Gamma(\beta_i)}\frac{\Gamma(\alpha_i+1-y_i)\Gamma(\beta_i+y_i)}{\Gamma(\alpha_i+\beta_i+1)}	& \label{equation:5}								
\end{align}

It follows from Equations \ref{equation:4} and \ref{equation:5} that

\begin{align}
	\label{equation:6}
	\frac{ \P(\theta=1 \given   \y)}{ \P(\theta=0 \given   \y)} > 1	&\iff\prod_{i=1}^{n}\frac{\Gamma(\alpha_i+y_i)\Gamma(\beta_i+1-y_i)}{\Gamma(\alpha_i+1-y_i)\Gamma(\beta_i+y_i)} > \frac{1-c}{c}	\notag \\
																																	&\iff \sum_{i=1}^{n} \log\left(\frac{\Gamma(\alpha_i+y_i)\Gamma(\beta_i+1-y_i)}{\Gamma(\alpha_i+1-y_i)\Gamma(\beta_i+y_i)}\right) > \log\left(\frac{1-c}{c}\right)
\end{align}

Observe that $g_i(y_i) := \log\left(\frac{\Gamma(\alpha_i+y_i)\Gamma(\beta_i+1-y_i)}{\Gamma(\alpha_i+1-y_i)\Gamma(\beta_i+y_i)}\right)$
is the contribution of expert $i$ to the decision rule.
This expression can be simplified: $g_{i}(0) = \log\left(\frac{\beta_i}{\alpha_i}\right)$ and $g_{i}(1) = \log\left(\frac{\alpha_i}{\beta_i}\right)$.
As a result, Equation \ref{equation:6} 
can be simplified: applying Notation \ref{notation:unknown} to Equation  \ref{equation:6}, one obtains
\begin{align}
 \label{equation:7}
 \frac{ \P(\theta=1 \given   \y)}{ \P(\theta=0 \given   \y)} > 1	&\iff \sum_{i=1}^{n} w_i^{**} y_i > C^{**}
\end{align}
Similarly,
\begin{align}
 \label{equation:8}
 \frac{ \P(\theta=1 \given   \y)}{ \P(\theta=0 \given   \y)} < 1	&\iff \sum_{i=1}^{n} w_i^{**} y_i < C^{**}
\end{align}

Theorem \ref{theorem:weighted_experts_unknown} follows from 
applying Equations \ref{equation:7} and \ref{equation:8} to 
Corollary \ref{Cor:bayes_flip} and Lemma \ref{Lemma:2}.

\subsection{Minimax decision rule}
\label{section:minimax_unknown}

Observe that the parameter space in this section is very conservative,
in the sense of including a wide variety of possibilities for the experts' accuracies.
For example, when every $\gamma_i = \frac{1}{2}$,
then the experts' opinions are distributed according to a $\mbox{Bernoulli}(\frac{1}{2})$,
that is, their opinions do not bring much information about $\theta$.
Therefore, it might not be surprising to find a minimax decision rule
that does not depend on the experts' opinions.
In this subsection we show that a coin flip, as in Definition \ref{definition:coin_flip},
is a minimax decision rule.

\begin{Definition}
 \label{definition:coin_flip}
 The coin flip is a decision rule, $\delta^{cf}$, such that, for every $\y$,
 \begin{align*}
  \delta^{cf}(\y, u) = u
 \end{align*}
\end{Definition}

\begin{thm}
 \label{theorem:minimax_coin}
 The coin flip, $\delta^{cf}$, is a minimax decision rule.
\end{thm}

The reason why a coin flip is, in the minimax sense,
a good decision is that, for every expert, $\Theta$
includes the possibility that the expert is accurate
and that he is inaccurate.
As a result, if the decison-maker decides
to use the experts' opinions, then there exists
a possibility that he is performing poorly.
The minimax perspective is conservative with this respect
and leads the decision-maker to avoid this undesirable possibility 
by ignoring the experts' opinions. 

In order to prove that the coin flip is a minimax decision rule,
we use Theorem \ref{theorem:minimax_technique} and apply the same strategy as 
when the experts' accuracies were known.
Assume that, for every $i$, $\alpha_i=\beta_i=1$ and $\P(\theta=1)=1/2$.
This assumption corresponds to a uniform distribution over $\Theta$. 
In this case, Equations \ref{equation:4} and \ref{equation:5} simplify to
\begin{align*}
 \begin{cases}
\P(\theta=1 \given \y)	&= \frac{1}{2\P(\Y=\y)}\prod_{i=1}^{n}\frac{\Gamma(1+y_i)\Gamma(2-y_i)}{2}	\\
\P(\theta=0 \given \y)	&= \frac{1}{2\P(\Y=\y)}\prod_{i=1}^{n}\frac{\Gamma(2-y_i)\Gamma(1+y_i)}{2}
 \end{cases}
\end{align*}

In other words, for every $\y$, 
$\P(\theta=1 \given \y) = \P(\theta=0 \given \y) = \frac{1}{2}$.
Therefore, when the distribution over $\Theta$ is uniform,
the (marginal) distribution \emph{a posteriori} for $\theta$ is the same
as the distribution \emph{a priori} for $\theta$.
Since the decision-maker is, \emph{a priori}, 
indifferent between $\theta=1$ and $\theta=0$ ($\P(\theta=0)=\P(\theta=1)$),
he is also indifferent between these options \emph{a posteriori},
no matter what the experts' opinions are.
It follows from Lemmas \ref{Lemma:1} and \ref{Lemma:2} 
that every decision rule is a Bayes decision rule.
In particular, according to Corollary \ref{Cor:bayes_flip}, the coin flip is a Bayes decision rule.
%The posterior distribution 
%depends on the data only through $\prod_{i=1}^{n}\Gamma(y_i+1)\Gamma(2-y_i)\Gamma^{-1}(3)$.
%As $\Gamma(y_i+1)\Gamma(2-y_i)\Gamma^{-1}(3)=1/2$ whether $y_i=1$ or $y_i=0$,
%the posterior distribution does not depend on the data. Therefore, for every $\y \in \{0,1\}^n$, $\P(\theta=1 \given  \y)=\P(\theta=1)=1/2,$
%that is, the data does not change one's uncertainty on $\theta$.
%It follows that  for any decision $a \in \{0,1\}$ and parameter value
%$\xi \in \Theta$, $\E[L(\xi,a) \given   \y]=1/2,$
%and hence both decisions are Bayes for any $\y \in \{0,1\}^n$.
%Hence, any decision rule $\delta(\y)$ is a Bayes decision against the uniform prior.
%In particular, the (randomized) rule that decides for $\theta=1$
%with probability half and for $\theta=0$ with probability half is a Bayes rule. 
%Moreover, its risk function
%is given by 
Also observe that, for every $\boldsymbol{\gamma}\in [0,1]^n$,
\begin{align*}
 \begin{cases} 
  R(\delta^{cf},(0,\boldsymbol{\gamma}))=\P(\delta^{cf}(\Y,U)=1 \given (0,\boldsymbol{\gamma})) = \P(U=1|\given (0,\boldsymbol{\gamma})) = \frac{1}{2} \\
  R(\delta^{cf},(1,\boldsymbol{\gamma}))=\P(\delta^{cf}(\Y,U)=0 \given (1,\boldsymbol{\gamma})) = \P(U=0|\given (1,\boldsymbol{\gamma})) = \frac{1}{2}
 \end{cases}
\end{align*}
Therefore, $\delta^{cf}$ is a Bayes decision rule with constant risk.
It follows from Theorem \ref{theorem:minimax_technique} that the coin flip is minimax.

The coin flip is not the unique minimax rule.
If a decision rule, $\delta^{*}$, is such that,
for every $\y$, $\delta^{*}(\y,U)$ is either $U$ or $1-U$,
then the risk function of $\delta^{*}$ is the same as that of $\delta^{cf}$.
As a result, such a $\delta^{*}$ is minimax.
Indeed, all of the minimax decision rules are 
of the form described above.
Let $\delta^{*}$ be a decision rule such that, for some $\y$,
$\delta^{*}(\y,U)$ is neither $U$ or $1-U$.
Since the image of $\delta^{*}$ is $\{0,1\}$,
then $\delta^{*}(\y,U)$ is some constant $d_{\y}$.
Let $\theta^{*} = 1-d_{\y}$.
Let $\gamma^{*}$ be such that, for every $i$, 
$\gamma_{i}^{*}=1$, if $y_{i}=\theta^{*}$,
and $\gamma_{i}^{*}=0$, otherwise.
Observe that, $P(\Y=\y|\gamma^{*}, \theta=\theta^{*})=1$,
that is, when $\theta=\theta^{*}$ and the expert precisions are $\gamma^{*}$,
the decision is $d_{\y} = 1-\theta^{*}$ and is incorrect with probability $1$.
As a result, the supremum of the risk of this decision rule is $1$ and 
it is not minimax.

The above derivation shows that the coin flip can also be justified in the Bayesian framework, 
when one adopts a uniform prior over $\Theta$,
the least favorable prior (definition \ref{definition:least_favorable} and theorem \ref{theorem:least_favorable}).
In this case, one not only is indifferent between $\theta=0$ and $\theta=1$ a priori but also assumes
that the experts' opinions bring no information about $\theta$, that is,
$P(\theta=1|\y)=\frac{1}{2}$, for every $\y$	.
This induces the coin flip to be as good as every other decision rule.
However, for other priors over $\Theta$,
the experts' opinions generally bring more information about $\theta$.
In such cases, using this information is better than the coin flip
and the latter is not a Bayes decision rule.
By including prior information about the accuracy of the experts,
the Bayesian decision rules are generally less conservative than the minimax decision rule.

\subsection{Frequentist \emph{a priori} information}

\begin{figure}[hp!]
 	\centering
 	\includegraphics[page=1,scale=0.4]{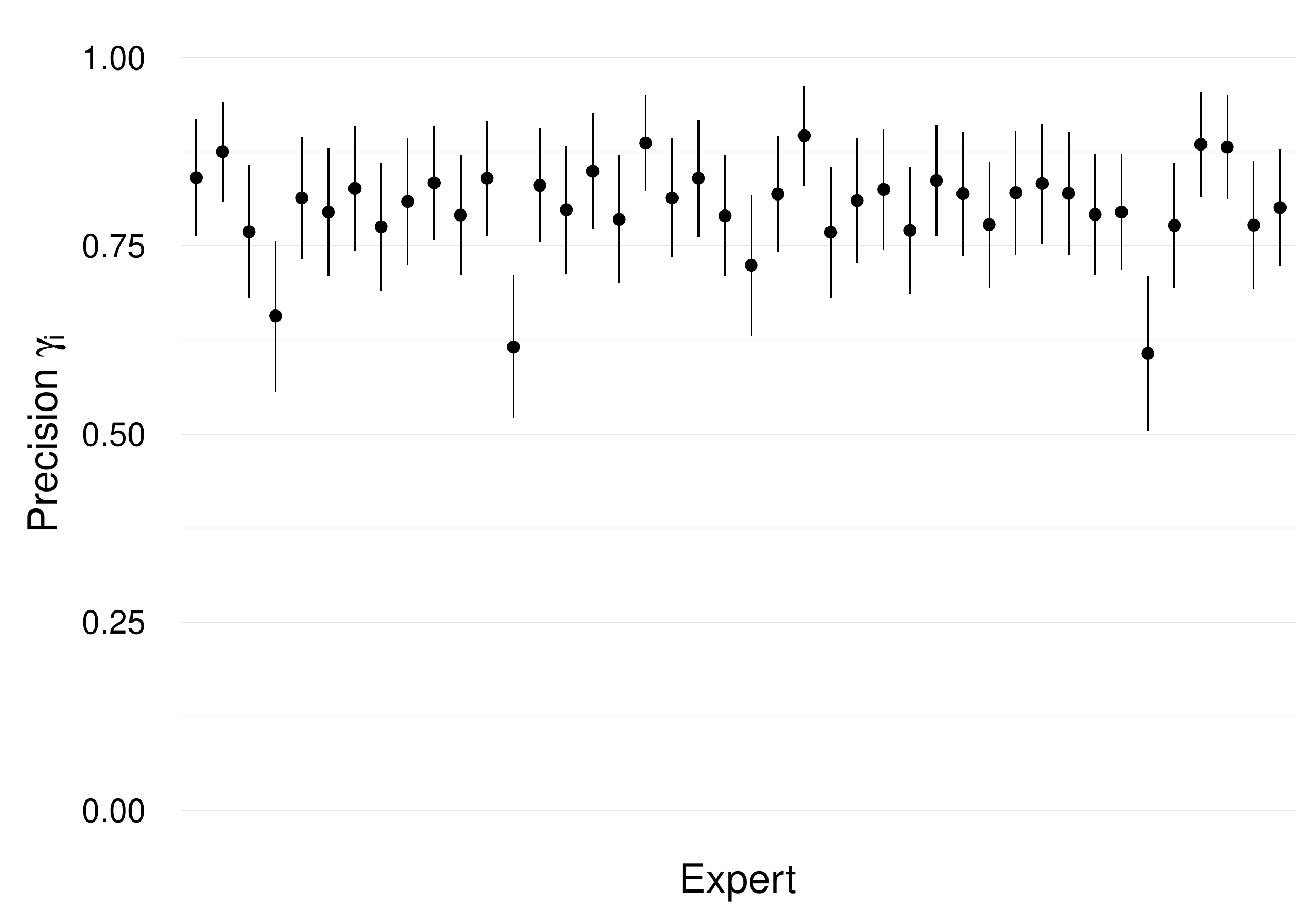} %\\
 	\caption{Estimated precision $\gamma_i$ of each of the 42 astronomers
 	from the  dataset analyzed by \citet{IzbickiStern}.}
 	\label{fig::astroEstimates}
 \end{figure}

Despite the conclusion in Theorem \ref{section:minimax_unknown}
that a coin flip is a minimax decision rule,
this decision rule is generally understood as a bad one.
One way to solve this apparent contradiction is to 
add \emph{a priori} information to the crowdsourcing problem by
restricting the parameter space.
Often, it is reasonable to assume that every expert is better than a coin flip
by at least some known threshold. For instance, for the cosmology
problem discussed in the introduction, the model from
\citet{IzbickiStern} provides the estimates for each of the  precisions $\gamma_i$ shown in
Figure \ref{fig::astroEstimates}, which are in agreement with such assumption.

Formally, this assumption can be written as
\begin{Assumption}
 \label{assumption:good_expert} 
 The parameter space of the crowdsourcing problem is
 \begin{align*}
  \Theta^*=\left\{(\theta,\gamma_1, \ldots, \gamma_n):\mbox{ }\theta \in \{0,1\}, \gamma_i \geq \frac{1}{2}+\epsilon \mbox{ }\forall i=1,\ldots,n\right\}
 \end{align*}
 for some \emph{fixed} $\epsilon$ such that $0<\epsilon<\frac{1}{2}$.
\end{Assumption}

\begin{thm}
 \label{thm:majority_unknown}
 Under Assumption \ref{assumption:good_expert},
 $\delta_{M}$ (corollary \ref{corollary:minimax_majority})
 is a minimax decision rule.
\end{thm}
That is, if every expert is better than a coin flip
by at least some known threshold, 
then the majority rule is once again minimax.

In order to prove Theorem \ref{thm:majority_unknown}, one cannot use the strategy from Section \ref{sec::minimaxFixed}.
This happens because, although $\delta_{M}$ is a Bayes decision rule, 
it doesn't have constant risk.
Indeed, the higher the accuracy of the experts, the smaller the risk of $\delta_{M}$.
In order to prove that $\delta_{M}$ is minimax,
we use the following standard strategy for finding minimax rules:
\begin{enumerate}
 \item Find a decision rule, $\delta_{M}$, that one believes to be minimax.
 \item Derive an upper bound on the minimax risk using $\delta_{M}$, that is: 
 $$\inf_{\delta} \sup_{\xi \in \Theta^*}\mbox{ }R(\delta,\xi)\leq\sup_{\xi \in \Theta^*}\mbox{ }R(\delta_{M},\xi)$$
 \item Derive a lower bound on the minimax risk by considering a subset of the parameter space.
 That is, if $\Theta^{*}_{0} \subset \Theta^{*}$, then
 $$\inf_{\delta} \sup_{\xi \in \Theta^*_0}\mbox{ }R(\delta,\xi)\leq \inf_{\delta} \sup_{\xi \in \Theta^*}\mbox{ }R(\delta,\xi)$$
\end{enumerate}

If the upper bound matches the lower bound, then $\delta_M$ is minimax:

\begin{thm}[\citep{tsybakov2008introduction}]
 \label{theorem:minimax_technique_2}
 If there exist $\delta_{M}$ and $\Theta^{*}_{0} \subset \Theta^{*}$ such that
 $$\sup_{\xi \in \Theta^*}\mbox{ }R(\delta_{M},\xi) \leq \inf_{\delta} \sup_{\xi \in \Theta^*_0}\mbox{ }R(\delta,\xi)$$
 then $\delta_{M}$ is a minimax decision rule.
\end{thm}

Since all experts are better than a coin flip on $\Theta^*$, 
a reasonable guess for a decision rule is
the majority rule, $\delta_{M}$,
that is defined in Corollary \ref{corollary:minimax_majority}.
Indeed, recall from Theorem \ref{theorem:minimax_fixed} and Corollary \ref{corollary:minimax_majority} that,
if the expert accuracies are known and
are all equally greater than $\frac{1}{2}$,
then the majority rule is a minimax decision rule.

Following the second step above, 
we bound the value of $\sup_{\xi \in \Theta^*}\mbox{ }R(\delta_{M},\xi)$. Observe that
\begin{align*}
&\sup_{\xi \in \Theta^*} R(\delta_{M},\xi) =																						&						\\
&\sup_{\xi \in \Theta^*} \{\P(\delta_{M}(\Y)=1 \given   \theta=0,\boldsymbol{\gamma})\I(\theta=0)+\P(\delta_{M}(\Y)=0 \given \theta=1,\boldsymbol{\gamma})\I(\theta=1) \} =							&						\\
&\sup_{\xi \in \Theta^*} \left\{ \P\left(\bar{\Y} > \frac{1}{2} \given[\Big]  \theta=0,\boldsymbol{\gamma}\right)\I(\theta=0)+\P\left(\bar{\Y} < \frac{1}{2} \given[\Big]  \theta=1,\boldsymbol{\gamma}\right)\I(\theta=1) \right\} =	&						\\
&\sup_{\xi \in \Theta^*}\P\left(\bar{\Y} > \frac{1}{2} \given[\Big] \theta=0, \boldsymbol{\gamma}\right) \leq														& \left[\bar{Y}|\theta=1 \stackrel{d}{=}	1-\bar{Y}|\theta=0\right] \\
&\P\left(\bar{\Y} > \frac{1}{2} \given[\Big]  \theta=0,\gamma_i=\frac{1}{2}+\epsilon\mbox{ }\forall i=1,\ldots,n\right).												&
\end{align*}

The last inequality follows from observing that 
$\P\left(\bar{\Y} > \frac{1}{2} \given[\Big] \theta=0, \boldsymbol{\gamma}\right)$
is strictly decreasing over each $\gamma_i$
(the higher the accuracies of the experts,
 the lower the value of $\bar{\Y}$ when $\theta=0$).
Therefore, this probability is maximized by picking, for each $i$,
the lowest value of $\gamma_i$ that is allowed on $\Theta^{*}$.
Using the above derivation, conclude that
\begin{align}
\label{eq::upper}
\sup_{\xi \in \Theta^*} R(\delta_{M},\xi) \leq \P\left(\bar{\Y} > \frac{1}{2} \given[\Big]  \theta=0,\gamma_i=\frac{1}{2}+\epsilon\mbox{ }\forall i=1,\ldots,n\right).
\end{align}

The third step in the derivation strategy obtains a lower bound on the minimax risk
by considering some subspace $\Theta_0^*$ of $\Theta^*$.
This subset is often called a \emph{reduction scheme} \citep{tsybakov2008introduction}.
Consider
$$\Theta_0^*=\left\{\left(0,\frac{1}{2}+\epsilon,\ldots,\frac{1}{2}+\epsilon\right),\left(1,\frac{1}{2}+\epsilon,\ldots,\frac{1}{2}+\epsilon\right)\right\}.$$

It follows from Theorem \ref{theorem:minimax_fixed} and Corollary \ref{corollary:minimax_majority}
that $\delta_{M}$ is a minimax decision rule 
over $\Theta^*_0$ and has constant risk.
Therefore,
\begin{align}
 \label{eq::lower}
 \inf_{\delta} \sup_{\xi \in \Theta_0^*}\mbox{ }R(\delta,\xi) =\sup_{\xi \in \Theta_0^*}\mbox{ }R(\delta_M,\xi)=\P\left(\bar{\Y} > \frac{1}{2} \given[\Big]  \theta=0,\gamma_i=\frac{1}{2}+\epsilon\mbox{ }\forall i=1,\ldots,n\right).
\end{align}

Finally, one obtains that $\sup_{\xi \in \Theta^*} R(\delta_{M},\xi) \leq \inf_{\delta} \sup_{\xi \in \Theta_0^*}\mbox{ }R(\delta,\xi)$ 
by combining Equations \ref{eq::upper} and \ref{eq::lower}.
Therefore, one concludes from Theorem \ref{theorem:minimax_technique_2} that
the majority rule, $\delta_{M}$, is 
a minimax decision rule 
considering $\Theta^*$ in place of the original $\Theta$.
Notice that a similar argument holds in the case where the 
accuracy threshold, $\epsilon$, varies among experts.
Let $\epsilon_i$ denote the accuracy threshold for expert $i$.
In this case, the decision rule $\delta^{*}_{\frac{1}{2}}$
from Theorem \ref{theorem:minimax_fixed} that uses
$\gamma_i = \frac{1}{2}+\epsilon_i$
is a minimax decision rule.
The minimax decision rule incorporates the information about the experts' precisions 
that is expressed in the choice of the parameter space.

\section{Conclusions}
\label{sec::concl}

We illustrate a set
of standard techniques for finding Bayes and minimax
decision rules in an elementary crowdsourcing example. 
The following techniques are illustrated:
\begin{itemize}
  \item \textbf{Theorem \ref{theorem:extensive_form}}:
  Rather than directly minimizing $\E_\theta[ R(\delta(\Y),\theta)]$,
  the Bayes rule can be found by minimizing $\E[L(a,\theta) \given \y]$ for each $\y$.
  \item \textbf{Theorem \ref{theorem:minimax_technique}}:
  Finding a Bayes rule with constant risk is a way
  to find a minimax rule, 
  and the prior associated to this
  rule is  the least favorable prior.
  \item \textbf{Theorem \ref{theorem:minimax_technique_2}}:
  A minimax rule can be found by providing an upper bound
  and a lower bound of the minimax risk and showing that they match.
  An upper bound on the minimax risk can be found by using the 
  supremum of the risk of a reasonable (candidate) decision rule.
  A lower bound on the minimax risk can be found using a reduction scheme.
\end{itemize}

The crowdsourcing problem was also used to
illustrate the differences between the Bayes decision rules
and the minimax decision rules.
Since a minimax decision rule considers the worst case scenario,
it is more conservative than a Bayes decision rule.
However, this conservatism can lead to odd minimax decision rules,
such as the coin flip in Theorem \ref{theorem:minimax_coin}.
This situation is overcome by adding more \emph{a priori}
information to the model.
Such information can be incorporated either as a
prior distribution over expert accuracies, in the Bayesian framework,
or as restrictions on the parameter space, in the frequentist framework.

The second author presented this material in a master's level decision theory class. 
The content was delivered in $2$ hours and 
the students were interested by the example, 
specially because of its practical appeal in comparison 
to the standard normal and binomial examples. 
They were particularly surprised by some situations then introduced.
For example, they were curious about the case in which minimax procedures reduce to flipping a coin
and about the situation in which every decision rule was a Bayes decision rule.

\bibliographystyle{plainnat}
\bibliography{paper}

\end{document}